\documentclass[12pt]{article}
\usepackage[dvips]{epsfig}
\usepackage{amsfonts}
\usepackage{amscd}
\usepackage{pstricks}
\font\tensym=msbm10
\font\sevensym=msbm7
\font\fivesym=msbm5

\font\tengoth=eufb10
\font\sevengoth=eufb7
\font\fivegoth=eufb5

\newfam\symfam
\textfont\symfam=\tensym
\scriptfont\symfam=\sevensym
\scriptscriptfont\symfam=\fivesym
\def\sym{\fam\symfam\tensym}

\newfam\gothfam
\textfont\gothfam=\tengoth
\scriptfont\gothfam=\sevengoth
\scriptscriptfont\gothfam=\fivegoth

\def\hs{\hbox to 3mm{}}
\def\hhs{\hbox to 5cm{}}
\def\ss{\smallskip}
\def\ms{\smallskip}

\def\JPicScale{0.8}\ifx\JPicScale\undefined\def\JPicScale{1}\fi

\def\A{\mathcal{A}}

\def\B{\mathbb{B}}

\def\N{\mathbb{N}}

\def\ga{\gamma}

\def\E{\mathcal{E}}
\def\HWS{{\cal H}_{WS}}
\def\1H{\mathbf{1}_{\HWS}}

\def\Q{{\sym Q}}

\def\al{\alpha}
\def\be{\beta}
\def\ep{\varepsilon}
\def\la{\lambda}

\def\om{\omega}

\def\2m#1#2{(#2 #1)}
\def\3m#1#2#3{(#3 #2 #1)}

\def\ncp#1#2{#1\langle #2\rangle}
\def\ncs#1#2{#1\langle\langle #2\rangle\rangle}
\newcommand{\End}{\mathrm{End}}
\def\scal#1#2{\langle #1 | #2 \rangle}
\def\ket#1{|#1\rangle}

\def\ra{\rightarrow}

\def\ol#1{\overline{#1}}

\def\adots{\mathinner{\mkern2mu\raise1pt\hbox{.}
\mkern3mu\raise4pt\hbox{.}\mkern1mu\raise7pt\hbox{.}}}

\def\up#1{\raise 1ex\hbox{\footnotesize#1}}

\def\mref#1{{\footnotesize ({\ref{#1}})}}

\newtheorem{expl}{Example}[section]

\newtheorem{note}[expl]{Note}

\newtheorem{theorem}[expl]{Theorem}

\newtheorem{definition}[expl]{Definition}
\newtheorem{proposition}[expl]{Proposition}

\newtheorem{remark}[expl]{Remark}

\newcommand\acknowledgements{\section*{Acknowledgements}}
\def\Proof{\medskip\noindent {\it Proof --- \ }}

\begingroup
\def\2#1{\ifnum#1<10 0\fi\the#1}
\xdef\isodayandtime{
{\2\day-\2\month-\the\year\space\2{\count0}:%
\2{\count2}}}
\endgroup

\title{\LARGE Sweedler's duals and Sch\"utzenberger's calculus}

\author{
\sc G\'erard H. E. Duchamp and Christophe Tollu.\rm
\thanks{LIPN - UMR 7030
CNRS - Universit\'e Paris 13
F-93430 Villetaneuse, France}}

\date{}
\begin{document}

\maketitle

\begin{abstract} We describe the problem of Sweedler's duals for bialgebras as essentially characterizing the domain of the transpose of the multiplication. 
This domain is the set of what could be called ``representative linear forms'' which are the elements of the algebraic dual which are also representative on 
the multiplicative semigroup of the algebra.\\ 
When the algebra is free, this notion is indeed equivalent to that of rational functions of automata theory. For the sake of applications, the range of 
coefficients has been considerably broadened, i.e. extended to semirings, so that the results could be specialized to the boolean and multiplicity cases. 
This requires some caution (use of ``positive formulas'', iteration replacing inversion, stable submodules replacing finite-rank families for instance). For the theory and its 
applications has been created a rational calculus which can, in return, be applied to harness Sweedler's duals. A new theorem of rational closure and 
application to Hopf algebras of use in Physics and Combinatorics is provided. The concrete use of this ``calculus'' is eventually illustrated on an example.
\end{abstract}

\section{Introduction}

This paper is entirely devoted to questions of rationality which arose, seemingly independently, in automata theory (Sch\"utzenberger's calculus) and in the dualization of multiplication (Sweedler's duals). As in the classical (univariate) case, rationality shows itself twofold : on functions and on expressions.\\ 
The rationality framework exposed here can be considered as the noncommutative analogue of linear recurrences \cite{Fl}. It is known, (see \cite{Co,La}), that it is equivalent to state that the coefficients of the Taylor expansion of a function satisfy a linear recurrence or that the generating function itself is rational ({\em i.e.} the quotient of two polynomials). This equivalence has a counterpart in the theory of rational expressions and this (nowadays classical) theory can be considered as ``localized at zero'' ({\em i.e.} analogous to the theory of rational functions
without a pole at zero). 

\ss
The paper is organized as follows.\\
In section 2, one sets out the theory of representative functions which were introduced in \cite{Sw} and are now standard in the theory of algebraic groups \cite{Ho}. Most of the material of this section (preparatory to the subsequent ones) is not new and can be found in several domains (but we believe that, beyond the needs of exposition, the description of the link itself will be of some use to these domains). The scope of proposition \mref{finite_orbit} can be harmlessly extended to semigroups \cite{Abe}.\\
As representative functions on the free monoid are the core of automata theory \cite{BR,Eil,S1,S2}, the domain of their scalars can (and, in fact, had to) be considerably enlarged to include structures allowing matrix (with unit) computations. This is the domain of semirings \cite{Go1,Go2}, the scalars of automata theory.\\ 
Section 4 prepares the link between representative functions and rational expressions by means of the notion of star (which is the ``positive analogue'' of the inverse \cite{Ca,He,Le,Pa,Ri}) and star closure (which is the analogue of rational closure in classical algebra\footnote{and ``sous-alg\`ebres pleines'' in spectral theories \cite{B_Th_Spec}.}).\\
The framework is then ready for a correct exposition of rational expressions which is the main concern of section 4.\\
In section 5, we apply what has been constructed to the dualization of bialgebras, then recovering known results.\\
In section 6, we show how to apply this ``rational calculus'' to solve the carrier problem in combinatorial physics.
 
\acknowledgements The authors wish to acknowledge support from Agence Nationale de la Recherche (Paris, France) under Program No. ANR-08-BLAN-0243-2 as well as support from ``Projet interne au LIPN 2009'' ``Polyz\^eta functions''.

\section{Representative functions on a semigroup}

The aim of this section is to discuss the dualization of bi-algebras and Hopf algebras. This problem, solved by Sweedler's duals, is the following.\\ 
Let $k$ be a field and $(B,.,\Delta,1_B,\ep)$ be a $k$-bialgebra ; we know that, if $B$ is finite-dimensional (resp. graded in finite dimensions), the dual (resp. graded dual) endowed with the transpose structure is a bialgebra and that, in case $B$ is a Hopf algebra this statement carries over. Now the question can be asked.\\
{\it What is the good notion of retricted dual for the general 
(i.e. ungraded finite or infinite dimensional) case ?}\\       
Analysing the dualization of the structural operations $(.,\Delta,1_B,\ep)$ of $B$, one sees at once that only the dualization of the multiplication is problematic as, in the general case, the codomain of the transpose of $.$ is larger than $B^*\otimes B^*$.\\ 
The first result follows (and somehow extends) \cite{Abe}. To state it, we need the notion of (left and right) shifts of functions on a semigroup. Let $k$ be a field, $(S,.)$ a semigroup and $f\in k^S$. For each $s\in S$ define 
$f_s\ :\ x\mapsto f(sx)$ (right shift of $f$) and $_sf\ :\ x\mapsto f(xs)$ (left shift of $f$), $_sf_t\ :\ x\mapsto f(txs)$ (bi-shift of $f$). Then, we have.
\newpage   

\begin{proposition}\label{finite_orbit} (see also \cite{Abe} section 2.2)
Let $k$ be a field, $(S,.)$ a semigroup and $f\in k^S$. The following are equivalent : \\
i) The family $(f_s)_{s\in S}$ is of finite rank in $k^S$\\
ii) The family $(_sf)_{s\in S}$ is of finite rank in $k^S$\\
iii) The family $(_sf_t)_{s,t\in S}$ is of finite rank in $k^S$\\
iv) There exists a double family $(g_i,h_i)_{1\leq i\leq n}$ of functions such that
\begin{equation}
\Big(\forall x,y\in S\Big)\Big(f(xy)=\sum_{i=1}^n\ g_i(x)h_i(y)\Big)
\end{equation}
v) There exists $\lambda\in k^{1\times n},\gamma\in k^{n\times 1}$ and $\mu\ :\ (S,.)\rightarrow (k^{n\times n},\times)$ a morphim of semigroups such that $(\forall s\in S)(f(s)=\lambda\mu(s)\gamma)$. 

\ss
Moreover, if $S$ admits a neutral (i.e. is a monoid), its image by $\mu$ of (v) above can be chosen to be the unity matrix.  
\end{proposition}
\Proof Omitted

\ms The elements of $k^S$ which fulfill the above conditions will be called {\it representative functions} on $S$ and denoted $R(k,S,.)$ \cite{Abe,CP,GOF21}.

\begin{remark} i) When $k$ is only a PID, proposition \mref{finite_orbit} above still holds with the five equivalent conditions, orbits and ranks being computed in $\bar k ^S$ ($\bar k$ is the fraction field of $k$).\\   
ii) If $S$ is finite, $R(k,S,.)=k^S$ and if $S$ is a group, one has
$$ 
R(k,S,.)=k^S\Longleftrightarrow S\textrm{ is finite}
$$
iii) If $S$ is a semigroup, the equivalences above are false in general as shown by the following counterexample. Let $G$ be a finite group and 
endow $S=\N\times G$ with the law $(n,g)*(m,h)=(0,gh)$. It can be easily checked that $(S,*)$ is a semigroup and that $R(k,S,*)=k^S$.\\ 
iv) When $W$ is a shift-invariant subspace of $k^S$ and $f\in W\cap R(k,S,.)$, the families 
$(f_s)_{s\in S}, (_sf)_{s\in S}, (_sf_t)_{s,t\in S}$ are of course in $W\cap R(k,S,.)$ (and are of finite rank). Two useful examples of such ``relative representative functional spaces'' are with $W=C(S)$ (continuous functions, $S$ being a topological semigroup) and $W=S^*$ (linear forms, $S$ being an algebra with its product as semigroup law).\\
v) If $T\subset S$ is a subsemigroup with finite set-theoretical complement ($S\backslash T$ is finite) then $f\in k^S$ is representative iff $f|_T$ is so. 
In particular, if $S=T^{(1)}$, the monoid obtained by adjunction of a unity to $T$, one has that $f\in k^S$ is representative iff $f|_T$ is.\\
vi) The proof of proposition \mref{finite_orbit} can also be found in \cite{Ho} where it does not use the structure of group.
\end{remark}

\section{Semirings}

Throughout the text ``monoid'' stands for ``semigroup with unit''.\\
Semirings are the structures adapted to matrix (with unity) computation. A semiring $(k,+,\times)$ consists of the following data 
\begin{itemize}
	\item[-] a set $k$
	\item[-] two binary laws $+,\times$ on $k$
\end{itemize}
such that
\begin{itemize}
	\item[-] $(k,+)$ and $(k,\times)$ are monoids, the first being commutative, their neutrals will be denoted respectively $0_k$ and $1_k$ 
	\item[-] $\times$ is left and right distributive over $+$
	\item[-] $0_k$ is an annihilator i. e. $(\forall x\in k)(0_k.x=x.0_k=0_k)$
\end{itemize}

\begin{expl} i) Any ring.\\
ii) The boolean semiring $\B=\{0,1\}$ endowed with the laws $x\oplus y=x+y-xy$ and $x\otimes y=xy$.\\ 
iii) The semiring $([-\infty,+\infty[,max,+)$, called in the literature ``(max,plus)-semiring''.\\
iv) In the semiring $([0,+\infty[,+,\times)$, the laws are continuous at infinity and then can be completed. We obtain a semiring $([0,+\infty],+,\times)$ which is suited for multiplicities arising in repeated additions of positive values during iterations. 
\end{expl}
\begin{expl} The following example is fundamental and will be used in the definition of CM-modules. Let $(M,+)$ be a commutative monoid, then $(\End(M),+,\circ)$ (defined as though $M$ were a group) is a semiring. The units are respectively, the constant mapping $M\ni m\mapsto 0_M$ for $+$ and $Id_M: m\mapsto m$ for $\circ$. 
\end{expl}

The structure of semiring defines a category larger than that of rings, the morphisms being defined similarly. Let $(k_i,+_i,\times_i),\ i=1,2$ be two semirings, a mapping $\phi : k_1\rightarrow k_2$ is called a morphism of semirings iff it is a morphim for the two structures of monoids (additive and multiplicative), then compatible with the laws and units of $k_1$ and $k_2$.

The definition of modules (here called CM-modules as they are constructed on Commutative Monoids as vector structure) follows also the classical pattern. 

The structure of a (left) $k$-CM-module is given by the following data
\begin{itemize}
	\item[-]  a commutative monoid $(M,+)$
	\item[-] a morphism (the scaling morphism) of semirings $s : k\rightarrow \End(M)$. 
\end{itemize}
The structure of (right) $k$-CM-module is defined by replacing $\End(M)$ by $\End^{op}(M)$ the opposite semiring (constructed with the opposite multiplicative law). Bi- and multimodules are defined as in \cite{B_Alg_III} and follow the general philosophy of ``structures with operators''. 

\begin{expl} Let $X$ be a set, then $k^X$, the set of all functions $X\rightarrow k$ is naturally endowed with a structure of $k-k$ bimodule defined as in the case when $k$ is a ring. So is $k^{(X)}$, the set of finitely supported functions of $k^S$ (actually a sub-$k$-$k$ bimodule of $k^S$). 
\end{expl}

The free monoid generated by a set $X$ (finite of infinite) is the set of words (i.e. finite sequences of elements of $X$ comprising the empty one denoted by $1_{X^*}$) endowed with the concatenation law.

\section{Shift operators and rational closure}

Let $(M,.)$ be a (commutative or not) monoid. 
For a function $f : M\rightarrow k$ and $a\in M$, we define the following shift operators \cite{Abe,BR,Eil,Ja}

\begin{itemize}
	\item[$\bullet$] $f_a : x\mapsto f(ax)$ (right shift)
	\item[$\bullet$] $ _af : x\mapsto f(xa)$ (left shift)	
\end{itemize}

We also have to describe the analog, for CM-algebras, of {\it full subalgebras} and {\it full subalgebraic closures} (see \cite{B_Th_Spec} Ch. 1.1.4) and this requires the notion of {\it summability} \cite{BR}.\\
\begin{definition}
	A family $(f_i)_{i\in I}$ of functions $M\rightarrow k$ is called {\rm summable} iff, for each $m\in M$, 
	$(f_i(m))_{i\in I}$ is finitely supported. Then, the mapping $m\mapsto \sum_{i\in I}f_i(m)$ is denoted 
	$\sum_{i\in I}f_i$ and called the sum of $(f_i)_{i\in I}$.  
\end{definition}

As a consequence, it is easily checked that, if $M$ is {\it locally finite} (\cite{Eil} Vol. A VII.4) and 
$f\ :\ M\rightarrow k$ is without constant term (i.e. $f(1_M)=0_k$), then the family $(f^n)_{n\in \N}$ (convolutional powers) is summable and its sum 
\begin{equation}
	\sum_{n\in \N}f^n
\end{equation}
will be denoted $f^*$ and called the {\it star} of $f$. 

\begin{note}
There is a lot of literature about the star problem (see \cite{Go1,Go2}). For a general discussion of star-type solutions in a semiring, see \cite{DHKL}.	
\end{note}

Now, we are in the position of stating the Kleene-Sch\"utzenberger theorem. 
  
\begin{theorem}\label{KS} Let $M=X^*$ be a free monoid, $k$ a semiring and $f\in k^M$. The following are equivalent\\ 
i) the family $(f_w)_{w\in M}$ belongs to a finitely generated shift-invariant left-submodule\\
ii) the family $(_wf)_{w\in M}$ belongs to a finitely generated shift-invariant right-submodule\\
iii) there exist a row $\la\in k^{1\times n}$, a column $\ga\in k^{n\times 1}$ and a representation (of monoids) $\mu\ :\ M\rightarrow (k^{n\times n},\times)$ such that $(\forall w\in M)(f(w)=\lambda\mu(w)\gamma)$.\\
\smallskip
If $X$ is finite, then (i-iii) above are also equivalent to:\\
iv) $f$ lies in the rational closure of $X$ (i.e. the smallest subalgebra of $k\langle\langle X \rangle\rangle$ closed under the star operation and containing $X$).
\end{theorem}
    
\begin{remark} i) Rational elements in the sense of (iii) infinitely many (linearly) independent shifts. That is why finitely generated shift-invariant submodules are needed in the general case. As an example one may consider 
$$
S=(a^*)^2\in \ncs{\N}{a}\ .
$$
One can check at once that $a^{-k}S=k.(a^*)+S$ and then for all $k$, $a^{-(k+1)}S\notin span(a^{-s}S)_{s=0}^k$.\\ 
ii) One can remove the hypothesis of freeness of $M$ if $k$ is a field. Indeed, in this case, the submodule can be taken as generated by the shifts (right or left) of $f$ and the representation is automatically compatible with the relations of $M$.\\
iii) Here the star is used as the localization at one (i.e. with positive formulas) of the inverse function. Indeed, with coefficients in a ring, if we are at the neighbourhood of $1$, the condition $(1-x)(1+y)=1$ (resp. $(1+y)(1-x)=1$) is equivalent to $y=x+xy$ (resp. $y=x+yx$). These self-reproducing positive conditions are taken as the definition of ``$y$ is a star of $x$'' in a semiring (see \cite{DHKL}).\\
iv) The condition (iv) in theorem \mref{KS} is known under the name of Kleene-Sch\"utzenberger theorem as, when $k$ is specialized to $\B$, this is actually Kleene's theorem. In this sense, this theorem lies at the frontier of harmonic analysis (the set of representative functions is dense in the Fourier space of compact groups), spectral theory (the notion of {\rm full subalgebra closure} comes from this theory \cite{B_Th_Spec}) and theoretical computer science (the notion of star was developped as a computational model of iteration and the notion of a semiring was developped to cope with general scalars as diverse as the ones arising in stochastic automata theory and shortest path problems).\\  
\end{remark}

In the general case ($X$ not necessarily finite), Kleene-Sch\"utzenberger's theorem has to be modified as follows.

\begin{theorem} Let $M=X^*$ be a free monoid, $k$ a semiring and $f\in k^M$. The following are equivalent\\ 
i) the family $(f_w)_{w\in M}$ belongs to a finitely generated shift-invariant left-submodule\\
ii) the family $(_wf)_{w\in M}$ belongs to a finitely generated shift-invariant right-submodule\\
iii) there exist a row $\la\in k^{1\times n}$, a column $\ga\in k^{n\times 1}$, and a representation (of monoids) $\mu\ :\ M\rightarrow (k^{n\times n},\times)$ such that $(\forall w\in M)(f(w)=\lambda\mu(w)\gamma)$.\\
iv) the function $f$ lies in the rational closure of $\ol{kX}=\{\sum_{x\in X}\al(x)x\}_{\al\in k^X}$ (i.e. the smallest subalgebra of 
$k\langle\langle X \rangle\rangle$ closed under the star operation which contains $\ol{kX}$).
\end{theorem}

\begin{remark} The rational closure of $X$ is, in fact, the intersection of the set of elements characterized by (i-iii) ({\em i.e.} Sweedler's dual of $\ncp{k}{X}$), and the algebra $\bigcup_{F\subset X\atop F\ finite}\ncs{k}{F}$ of the series whose support involves a finite alphabet.
\end{remark}

\section{Rational expressions}

The construction of \cite{CD} was localized at zero, we extend it here to any localization i.e. for any mapping $\Lambda\ :\ X\rightarrow k$.\\ 
As the rational closure involves a unary law (the star) partially defined, the definition of universal formulas for this closure needs some caution. Indeed, we need to build in parallel a ``character'' (the constant term) $const$ so that all proper expressions should have a star.

\ss
One first defines, as in \cite{CD} the completely free expressions (or formulas) as the terms of the universal algebra defined on $X\cup \{0_E\}$ ($0_E$, 
which does not belong to $X$ will serve as a null or void expression and will be mapped to the zero series). This algebra will be denoted $\E^{cf}(X,k)$. 
More precisely 

\ss
\begin{itemize}
\item[-] If $x\in X\cup\{0_\E\}$ then $x\in \E^{cf}(X,k)$.
\item[-] If $E,E_1,E_2\in \E^{cf}(X,k),$ and $\la \in k$ then
$$
\begin{array}{c}
E_1+E_2\in \E^{cf}(X,k),\
E_1 \cdot E_2\in \E^{cf}(X,k)\\
\la E\in \E^{cf}(X,k),
E\la\in \E^{cf}(X,k).\\
\end{array}
$$
\item[-] If $E\in \E^{cf}(X,k)$ then $E^*\in \E^{cf}(X,k)$.
\end{itemize}

The partial function $const\ :\ \E^{cf}(X,k)\rightarrow k$ (constant term) is constructed as follows: 

\begin{enumerate}
	\item If $x\in X$ $const(x)=\Lambda(x)$ and $const(0_\E)=0_k$.
	\item If $E,E_i\in \E^{cf}(X,k),\ i=1,2$ and $\la \in k$ then
$$
\begin{array}{c}
const(E_1+E_2)=const(E_1)+const(E_2),\
const(E_1 \cdot E_2)=const(E_1) \cdot const(E_2)\\
const(\la E)=\la const(E),\ 
const(E\la)=const(E)\la.\\
\end{array}
$$ 
\item If $const(E)=0_k$ then $const(E^*)=1_k$.
\end{enumerate}

The domain of $const$ will be called {\it rational expressions} and denoted $\E_\Lambda(X,k)$. For example 
$0_\E^*\in \E_\Lambda(X,k)$.

Let now $\Theta\ :\ X\rightarrow \ncs{k}{A}$ be a mapping such that, for every $x\in X$, the constant term of $\Theta(x)$, {\em i.e.} the coefficient of $1_{A^*}$ in $\Theta(x)$, is equal to $\Lambda(x)$ (in symbols, $\forall x\in X), [1_{A^*}]\Theta(x)=\Lambda(x)$). 
Following recursively (1-2-3) above, we can construct a polymorphism $\phi_{\Theta}\ :\ \E_\Lambda(X,k)\rightarrow \ncs{k}{A}$ which is a morphism for the laws (2 internal and 2 external) and the star. Moreover 
$\delta_{1_{A^*}}\circ \phi_{\Theta}=const$ (i.e. $const$ can be considered as a ``constant term function'' for the expressions). The image of $\phi_{\Theta}$ is exactly the rational closure of the set $\{\Theta(x)\}_{x\in X}$.   

\section{Dual laws and bialgebras}

Let $\Delta\ :\ \A\rightarrow \A\otimes \A$ be any comultiplication (i.e. $\A$ is a $k$-coalgebra). It is known that its dual $(\A^*,^t\Delta)$ is an algebra and if $\A$ is coassociative (resp. cocomutative, counital), $\A$ is associative (resp. commutative, unital) \cite{Abe}.\\
We would like here to enlarge the framework of \cite{DFLL}.\\ 
If $\A$ is an algebra, let us call \textit{dual law} on $\A^*$ a law of the form $^t\Delta$ for some (not necessarily coassociative) comultiplication on $\A$.\\ 
In \cite{DFLL} were considered the \textit{dual laws} on $\ncs{k}{X}\simeq\ncp{k}{X}^*$ in order to prove that the Hadamard and Infiltration products, which were known to preserve rationality, were essentially the only (along with an interpolation between the two) alphabetic (associative and unital) dual laws between series. The notion of \textit{dual law} provides an  implementation scheme for the automata so that the rationality preservation is naturally effective.

\begin{theorem} Let $\A$ be a $k$-algebra and $\Delta\ :\ \A\rightarrow \A\otimes \A$ be a comultiplication which is a morphism of algebras. Then \\
i) If $k$ is a field, Sweedler's dual $\A^\circ$ of $\A$ is closed under the dual law $^t\Delta$.\\
ii) If $k$ is a semiring and $\A=\ncp{k}{X}$, $\A^\circ$ is closed under the dual law $^t\Delta$.    
\end{theorem}
 
\begin{note} i) The theorem is no longer true if $\Delta\ :\ \A\rightarrow \A\otimes \A$ is arbitrary (i.e. not necessarily a morphism) as shows the following counterexample. With $\Delta\ :\ \Q[x]\rightarrow \Q[x]\otimes \Q[x]$ such that 
$\Delta(x)=\frac{1}{n!}(x^n\otimes x^n)$, one has 
$$
^t\Delta(\frac{1}{1-x},\frac{1}{1-x})=exp(x).
$$  
ii) In (i) above, the restriction on scalars (to be a field) can be extended to inductive limits of PIDs.\\
iii) Comultiplications different from morphisms can preserve rationality. For example, let $\Delta\ :\ \ncp{k}{X}\rightarrow \ncp{k}{X}\otimes \ncp{k}{X}$ be a morphism and $\Delta_1\ :\ \ncp{k}{X}\rightarrow \ncp{k}{X}\otimes \ncp{k}{X}$ be a linear mapping which coincides with $\Delta$ except for a finite number of words of $X^\star$. It can be checked that $\Delta_1$, although not a morphism, preserves rationality.
\end{note}


\medskip
Let us now return to the case of a bialgebra $(B,.,\Delta,1_B,\ep)$. The following proposition says that, if a linear form on $B$ is transformed by $^tm$ ($m : B\otimes B\rightarrow B$ is just the multiplication mapping) into an element of $B^*\otimes B^*$, it must be of the form exhibited in \textit{(iv)} of Proposition \ref{finite_orbit}. Let us state that precisely.\\

\begin{proposition} Let $(B,m,\Delta,1_B,\ep)$ be a bialgebra and $f\in B^*$.\\ 
1) The following are equivalent\\
i) $^tm(f)\in B^*\otimes B^*$ (for the canonical embedding $B^*\otimes B^*\hookrightarrow (B\otimes B)^*$)\\
ii) $f\in R(k,B,.)$\\
iii) $ker(f)$ contains a finite-codimension one-sided ideal\\
iv) $ker(f)$ contains a finite-codimension two-sided ideal\\
v) There exist $\lambda\in k^{1\times n},\gamma\in k^{n\times 1}$ and 
$\mu\ :\ (B,+,.)\rightarrow (k^{n\times n},+,\times)$ a morphim of $k$-algebras (associative with units) such that $(\forall x\in B)(f(x)=\lambda\mu(x)\gamma)$.\\
2) Moreover, let $B^0$ be the set of linear forms which are decomposable as $1-(i)$ above. Then $(B^0,^t\Delta,^tm,^t\ep,^t1_B)$ is a bialgebra, and if $B$ admits an antipode $\sigma$ (i.e. is a Hopf algebra), one has $^t\sigma(B^0)\subset B^0$ and $(B^0,^t\Delta,^tm,^t\ep,^t1_B,^t\sigma)$ is a Hopf algebra.   
\end{proposition}   

\section{An application of rational expressions to Combinatorial Physics}

In a joint work with J. Katriel, one of us gave the solution of
the problem of matrix coefficients in the Fock space of carriers between
two levels. We give here a brief review of these continued fractions-type formulas and provide a sketch of their proof.

In \cite{D20} was considered, as a Fock space, a general vector space $V$ over a field $k$ with basis 
$|e_n\rangle \;\; n= 0,\, 1,\, \cdots$, equipped with its natural grading
\begin{equation}
V=\oplus_{n\in {\bf Z}}V_n\ {\rm with}\ V_n:=ke_n;
\ V_{-n-1}:=\{0\}\ {\rm for}\ n\geq 0
\end{equation}
and scalar product defined by $\scal{e_n}{e_m}=\delta_{n,m}$.
Let $f,\ g$ be two linear operators on $V$ of degrees $-1,\ +1$, 
respectively. Generically, they read
\begin{equation}
\ {\rm for}\ n\geq 0;\ f|e_0\rangle :=
0;\ f|e_{k+1}\rangle =\alpha_{k+1}|e_k\rangle;
\ g|e_k\rangle :=\beta_{k+1}|e_{k+1}\rangle 	
\end{equation}
We consider the words in $f,\ g$:
\begin{equation}
w(f,g):=f^{p_1}g^{q_1}f^{p_2}g^{q_2}\cdots f^{p_n}g^{q_n}	
\end{equation}
the degree (excess) $\pi_e(w)$ of which is $\sum_{k=1}^n (q_k-p_k)$ (this is, for algebraists, the degree of the graded operator $f^{p_1}g^{q_1}f^{p_2}g^{q_2}\cdots f^{p_n}g^{q_n}$). 
This provides a representation $\mu$ of a two-letter free monoid $\{b_-,b_+\}^*$ on $V$ by 
$\mu(b_-)=f;\ \mu(b_+)=g$, which is graded for the weight on $\{b_-,b_+\}^*$. In order to keep the reading of a word 
from left to right, one performs the action on the right. Thus $V$ becomes a $\{b_-,b_+\}^*$ right module by 
\begin{equation}
e_0.b_-=0\ ;\ e_{n+1}.b_-=\al_{n+1}e_n\ ;\ e_{n}.b_+=\be_{n+1}e_{n+1}, \,\, n\geq 0
\end{equation}
one is interested by the matrix elements\footnote{An interpretation of a similar coefficient in terms of paths, namely for the computation of 
$$
\int_{-\infty}^{+\infty} x^ip_np_mw(x)dx
$$
where $(p_n)_{n\in\N}$ is a family of orthogonal polynomials for the weight $w(x)$ can be found in \cite{Viennot}.}.
\begin{equation}
	\scal{e_n.\{b_-,b_+\}^i}{e_m}=\om_{n\ra m}^{(i)}\ .
\end{equation}

Define $W_{n\ra m}^{(i)}$ as follows 
\begin{equation}
	W_{n\ra m}^{(i)}=\{w\in \{b_-,b_+\}^i\ |\ (\pi_e(w)=m-n)\ \mathrm{and}\ (w=uv\Longrightarrow \pi_e(u)\geq -n)\}\ .
\end{equation}

The following proposition characterizes $W_{n\ra m}^{(i)}$ as an universal transporter between level $n$ and level $m$.
 
\begin{proposition} 
i) If all the weights $\al_n;\ n\geq 1$, $\be_n, n\geq 0$ are nonzero, $W_{n\ra m}^{(i)}$ is exactly the set of 
words of length $i$ such that $\scal{e_n.w}{e_m}\not=0$.\\
ii) In all cases the latter is a subset of $W_{n\ra m}^{(i)}$ i. e. 
\begin{equation}
\{w\in \{b_-,b_+\}^i\ |\ \scal{e_n.w}{e_m}\not=0\}\subset W_{n\ra m}^{(i)} 	
\end{equation}
\end{proposition}

Indeed, $W_{n\ra m}^{(i)}$ admits factorizations as ``noncommutative continued fractions''. In order to do this, we recall for the reader the definition 
of Dyck, positive Dyck and negative restricted (by depth) Dyck codes, respectively: 
\begin{eqnarray}
D=\{w\in \pi_e^{-1}(0)| w=uv\textrm{ and } \ep\notin \{u,v\}\Longrightarrow \pi_e(u)\not=0\}\cr
D_+=\{w\in D | w=uv\textrm{ and } \ep\notin \{u,v\}\Longrightarrow \pi_e(u)> 0\}\cr
D_-^{(n)}=\{w\in D_-\ |\ \min_{uv=w}(\pi_e(u))=-n\}\ .
\end{eqnarray} 

Now, identifying the subsets defined previously with their characteristic series in $\ncs{\N}{b_-,b_+}$, one has the following factorizations (here,  equalities between rational expressions of $b\pm,D_+$ and $D_-^{(n)}$).\\
\begin{proposition}\label{fact_cont_frac} Set $W_{n\ra m}=\sum_{i\geq 0}W_{n\ra m}^{(i)}$.\\
1) One has the following factorizations for $W_{n\ra n+k}$
\begin{eqnarray}
(D_++D_-^{(n)})^*(b_+D_+^*)^k&=&\Big[\stackrel{\ra}{\prod}_{i=0}^{(k-1)}((D_-^{(n+i)})^*b_+)\Big](D_++D_-^{(n+k)})^*\textrm{ if } k\geq 0\cr
(D_+b_-)^k(D_++D_-^{(n+k)})^*&=&(D_++D_-^{(n)})^*\Big[\stackrel{\ra}{\prod}_{i=1}^{(-k)}((b_-(D_-^{(n-i)})^*)\Big] \textrm{ otherwise } 
\end{eqnarray}
2) One has the following self-reproducing equations for the Dyck codes
\begin{equation}
	D_+=b_+(D_+)^*b_-\ ;\ D_-^{(n)}=b_-(D_-^{(n-1)})^*b_+\ ;\ D_-^{(0)}=\emptyset
\end{equation}
\end{proposition}
The statements of this  proposition can be considered as a ``noncommutative continued fraction'' expansion of $W_{n\ra m}$.\\
As a result \cite{D20}, with 
\begin{equation}
\mu(W_{m\rightarrow n}^{(i)})\ket{n}=\omega_{n\ra m}^{(i)}\ket{m}\ \textrm{and}\ 
T_{n\rightarrow n+k}:=\sum_{i\geq 0}t^i\omega_{n\rightarrow n+k}^{(i)}\ ,
\end{equation}
using proposition \mref{fact_cont_frac} and the representations $\mu_t(b\pm)=tb\pm$ (observe that $\mu_1=\mu$), one can expand $T_{n\rightarrow n+k}$ as a product of continued fractions. Let 
\begin{eqnarray}\label{eq9}
                         F^+_n= {1\over\displaystyle 1-
{\strut t^2\alpha_{n+1}\beta_{n+1}\over\displaystyle 1-
{\strut t^2\alpha_{n+2}\beta_{n+2}\over\displaystyle 1-
{\strut t^2\alpha_{n+3}\beta_{n+3}\over\displaystyle 1-\cdots
}}}}={1\over 1-E^+_n}\cr\cr\cr
                         F^-_n= {1\over\displaystyle 1-
{\strut t^2\alpha_{n  }\beta_{n  }\over\displaystyle 1-
{\strut t^2\alpha_{n-1}\beta_{n-1}\over\displaystyle 1-
{\strut t^2\alpha_{n-2}\beta_{n-2}\over\displaystyle 1-\cdots
}}}}={1\over 1-E^-_n}
\end{eqnarray}
and 
\begin{equation}\label{eq10}
F_n= {1\over 1-E^+_n-E^-_n}\ .	
\end{equation}
Then, if $k\geq 0$, we obtain
\begin{equation}\label{eq11} 
T_{n\rightarrow n+k}=t^kF_{n+k}\prod_{i=0}^{k-1}F^-_{n+i}
=t^kF_n\prod_{i=1}^kF^+_{n+i}
\end{equation}
and, if $k\leq 0$
\begin{equation}\label{eq12}
T_{n\rightarrow n+k}=t^{-k}F_{n+k}\prod_{i=0}^{-k-1}F^+_{n-i}
=t^{-k}F_n\prod_{i=1}^{-k}F^-_{n-i}
\end{equation}

\section{Conclusion}

The dualization problem solved by Sweedler's duals has striking relations with language theory (for this last point, see \cite{BR}). 
A true ``engineer-like'' calculus was developped in order to handle the rational closure mentionned above. This set of formulas is mainly based on a recursion to compute the ``star of a matrix'' (the formulas and a complete discussion can be found in \cite{DHKL} and are reminiscent of general formulas giving the inverse of a matrix decomposed in blocks). This calculus is powerful enough to be the main ingredient in investigating rationality properties within various domains (see \cite{Co,D21} for noncommutative geometry, \cite{Cau} for functions on the free group and \cite{BR,Eil} for automata theory) and sufficiently expressive to give exact developments of some transfer coefficients in Combinatorial Physics \cite{D20}. These rational expressions are generic in the sense that any Sweedler's dual can be described by them. In this way, as for Dirac's notation, we can hope to inherit the computational skill developped through fourty years of practice.

\newpage


\begin{thebibliography}{ABC}
%
\bibitem{Abe} E. Abe, {\it Hopf algebras}. Cambridge Univ. Press, 1980.
%
\bibitem{BR} J. Berstel, C. Reutenauer, {\it Rational series and their 
languages}. EATCS Monographs on Theoretical Computer Science, Springer, 1988.
%
\bibitem{B_Alg_III} N. Bourbaki, {\it Algebra, chapter III}, Springer (1970)
%
\bibitem{B_Th_Spec} N. Bourbaki, {\it Th\'eories spectrales}, Hermann (1967)
%
\bibitem{Cau} G. Cauchon, {\it S\'eries de Malcev-Neumann sur le groupe 
libre et questions de rationalit\'e}, Theoret. Comp. Sci. {\bf 98} (1992) 79-97.
%
\bibitem{Ca} A. Cayley, {\it On certain results related to quaternions}, Phil. Mag. {\bf 26} (1845) 141-145.
%
\bibitem{CD} J.-M. Champarnaud, G. Duchamp, {\it Derivatives of rational expressions \-and \-re\-lated theorems},  Theoret. Comp. Sci. {\bf 313} (2004) 31.
%
\bibitem{CP} V. Chari, A. Pressley, {\it A guide to quantum groups.} Cambridge Univ. Press, 1994.
%
\bibitem{Co} A. Connes, {\it Noncommutative geometry.} Acad. Press, 1994.
%
\bibitem{D21} G. Duchamp, C. Reutenauer, {\it Un crit\`ere de rationalit\'e 
provenant de la g\'eom\'etrie noncommutative,} Inventiones Mathematicae, {\bf 128} (1997) 613-622.
%
\bibitem{DHKL} G. Duchamp, H. H. Kacem, \'E. Laugerotte, 
{\it Algebraic elimination of $\epsilon$-transitions}, DMTCS, \textbf{7} (2005) 51-70.
%
\bibitem{DFLL} G. Duchamp, M. Flouret, \'E. Laugerotte, J.-G. Luque, {\it Direct and dual laws for automata with multiplicities}, Theoret. Comp. Sci. {\bf 267} (2001) 105-120.
%
\bibitem{GOF21} G. H. E. Duchamp, P. Blasiak, A. Horzela, K. A. Penson, A. I. Solomon, 
{\it Hopf Algebras in General and in Combinatorial Physics: a practical introduction},\\ 
arXiv : {\tt 0802.0249}
%
\bibitem{Eil} S. Eilenberg, {\it Automata, languages and machines}. Acad. Press, New-York, 1974.
%
\bibitem{F1} M. Fliess, Sur le plongement de l'alg\`ebre des s\'eries 
rationnelles non commutatives dans un corps gauche. {\it CRAS Ser. A} {\bf 271} 
 (1970) 926-927.
%
\bibitem{Fl} M. Fliess, Matrices de Hankel. {\it Jour. of Pure and Appl. 
Math.} {\bf 53} (1994) 197-222.
%
\bibitem{Go1} Golan J. S., {\it Power Algebras over Semirings with Applications in Mathematics and 
Computer science}. Kluwer Academic Publishers, 1999.
%
\bibitem{Go2} Golan J. S., {\it Semirings and Affine Equations over Them: Theory and Applications}.  
Kluwer Academic Publishers, 2003.
%
\bibitem{He} A. Heyting, Die Theorie der linearen Gleichungen in einer 
Zahlenspezies mit nichtkommutativer Multiplikation. {\it Math. Ann.} {\bf 98} (1927) 
465-490.
%
\bibitem{Ho} G. P. Hochschild, Basic theory of algebraic groups and Lie algebras,
Springer 1981,
%
\bibitem{Ja} G. Jacob, {\it Repr\'esentations et substitutions matricielles 
 dans la th\'eorie matricielle des semigroupes}, Th\`ese, Univ. de Paris (1975).
%
\bibitem{D20} J. Katriel, G. Duchamp, {\it Ordering relations for q-boson operators, continued fractions techniques, 
and the q-CBH enigma}, J. Phys. A: Math. Gen. {\bf 28} (1995) 7209-7225.
%
\bibitem{La} S. K. Lando, {\it Lectures on generating functions}, A. M. S. (2003).
%
\bibitem{Le} J. Lewin, {\it Fields of fractions for group algebras of free groups},
 {\it Trans. Amer. Math. Soc.} {\bf 192} (1974) 339-346.
%
\bibitem{Pa} D.S. Passman, {\it The algebraic structure of group rings},  
John Wiley - Interscience, (1977). 
%
\bibitem{Ri} A.R. Richardson, {\it Simultaneous linear equations over a division ring}, Proc. Lond. Math. Soc. {\bf 28} (1928) 395-420. 
%
\bibitem{S1} M.P. Sch\"utzenberger, {\it On the definition of a family of 
automata}, Information and Control, {\bf 4} (1961) 275-270.
%
\bibitem{S2} M.P. Sch\"utzenberger, {\it On a theorem of R. Jungen}, Proc. Amer. Math. Soc. {\bf 13} (1962) 885-889.
%
\bibitem{Sw} M.E. Sweedler, {\it Hopf algebras}, W.A. Benjamin, New York, 1969.
%
\bibitem{Viennot} X.G. Viennot, {\it Une th\'eorie combinatoire des polyn\^omes orthogonaux}, Lect. Notes LACIM UQAM, Montreal (1984).\\
{\tt http://web.mac.com/xgviennot/iWeb/Xavier\_Viennot}
\end{thebibliography}
\end{document}